\documentclass[12pt]{amsart}

\def\@typesizes{%
       \or{5}{6.5}\or{6}{7.5}\or{7}{8.5}\or{8}{11}\or{9}{12}%
       \or{10}{13}
       \or{\@xipt}{14}\or{\@xiipt}{15}\or{\@xivpt}{18}%
       \or{\@xviipt}{20}\or{\@xxpt}{24}}

\usepackage{amsthm}
\usepackage{amsmath}
\usepackage{amssymb}
\usepackage{graphicx}
\usepackage{enumerate}
\usepackage{color}
\allowdisplaybreaks
\numberwithin{equation}{section}
\numberwithin{figure}{section}

\theoremstyle{plain}
\newtheorem{theorem}{ Theorem}[section]
\newtheorem{proposition}[theorem]{ Proposition}
\newtheorem{lemma}[theorem]{ Lemma}
\newtheorem{corollary}[theorem]{ Corollary}
\newtheorem{example}[theorem]{ Example}
\newtheorem{remark}[theorem]{ Remark}
\newtheorem{definition}[theorem]{ Definition}
\newtheorem{conjecture}{ Conjecture}

\def\BET{\begin{theorem}}
\def\ENT{\end{theorem}}
\def\BEP{\begin{proposition}}
\def\ENP{\end{proposition}}
\def\BEL{\begin{lemma}}
\def\ENL{\end{lemma}}
\def\BEC{\begin{corollary}}
\def\ENC{\end{corollary}}
\def\BEE{\begin{example} \rm}
\def\ENE{\end{example}}
\def\BER{\begin{remark} \rm}
\def\ENR{\end{remark}}
\def\BED{\begin{definition} \rm}
\def\END{\end{definition}}
\def\BECJ{\begin{conjecture}}
\def\ENCJ{\end{conjecture}}

\def\bea{\begin{eqnarray}}
\def\eea{\end{eqnarray}}

\def\beas{\begin{eqnarray*}}
\def\eeas{\end{eqnarray*}}

\def\beq{\begin{equation}}
\def\eeq{\end{equation}}

\def\beal{\begin{align*}}

\def\eeal{ \end{align*} }

\def\row{\nonumber \\ & & }
\def\roweq{\nonumber \\ &=& }

\def\rowpl{\nonumber \\ & \ \ + & }

\def\bff{{\bf f}}

\def\bfu{{\bf u}}

\def\bfF{{\bf F}}
\def\bfG{{\bf G}}
\def\bfH{{\bf H}}

\def\bfL{{\bf L}}

\def\bfX{{\bf X}}
\def\bfY{{\bf Y}}

\def\bbC{{\mathbb C}}

\def\bbR{{\mathbb R}}

\def\cA{{\mathcal A}}

\def\cL{{\mathcal L}}

\def\cP{{\mathcal P}}

\def\cR{{\mathcal R}}

\def\cX{{\mathcal X}}
\def\cY{{\mathcal Y}}

\def\ef{\eqref}

\begin{document}

\title[Differential Operators with Periodic Coefficients]{Floquet Problem and Center Manifold Reduction for Ordinary Differential Operators
with Periodic Coefficients in Hilbert Spaces}

\author{Vladimir Kozlov$^1$, Jari Taskinen$^2$}

\begin{abstract}
A first order differential equation with a periodic operator
coefficient acting in a pair of Hilbert spaces is
considered. This setting models both elliptic equations with
periodic coefficients in a cylinder and parabolic equations
with time periodic coefficients. Our main results are a
construction of a pointwise projector and
a spectral splitting of the system into a finite
dimensional system of ordinary differential equations with
constant coefficients and an infinite
dimensional part whose solutions have better properties in a
certain sense. This  complements
the well-known asymptotic results for periodic
hypoelliptic problems in cylinders \cite{KuchBook} and for elliptic problems in quasicylinders \cite{na417}.

As an application we give a center
manifold reduction for a class of non-linear ordinary differential
equations in Hilbert spaces with periodic coefficients. This result
generalizes the known case with constant coefficients, \cite{M1}, \cite{M}.
\vspace{2mm}

\end{abstract}

\maketitle

\begin{center}
$^1${\it Department of Mathematics, Link\"oping University, S--581 83 Link\"oping,
Sweden \\ $^2$ Department of Mathematics and Statistics, University of Helsinki, P.O.Box 68, 00014 Helsinki,
Finland} \\

\vspace{2mm}

E-mail: vladimir.kozlov@liu.se\,/\,jari.taskinen@helsinki.fi
\end{center}

\section{Introduction} \label{sec1}

Consider a first order ordinary differential equation (ODE) for an unknown function $x(t)$ with values in an infinite dimensional Hilbert space $Y$,
\begin{equation}\label{I1aa}
\frac{d x (t) }{dt}=Ax(t) +f(t;x(t)) , \ t\in\bbR,
\end{equation}
where $A$ is an  unbounded linear operator in $Y$, which is constant in $t$,
and  $f : \bbR\times D_A \to Y$ is given. If $\cP$ is a finite dimensional
orthogonal projector in  $Y$ which commutes with $A$, then the system
(\ref{I1aa}) with $f \equiv 0$ can be split into a finite dimensional system
on the  subspace $\cP(Y)$   and an infinite dimensional system which may have
better properties than the initial one. This reduction can be quite useful
in the study of the large time behaviour of linear dynamical systems
perturbed by a  linear or non-linear perturbation $f$ (see \cite{KM1}, \cite{M} and references there). The main
subject of this paper is to study similar splitting for the case when
$A=A(t)$ is a periodic operator function with certain Fredholm properties.
Our goal is to define
a pointwise projector $\cP=\cP(t)$, which is a projector for any $t \in \bbR$ and  commutes with the operator of the periodic
problem. 
This projector leads to the reduction of
the problem to a finite dimensional problem with time independent
operator and infinite dimensional problem having better properties than the original one. Another goal is to present a center manifold reduction for nonlinear periodic operators.

We recall the classical Floquet theorem
concerning  a system of ODEs
\begin{equation}\label{I1}
\frac{d x (t) }{dt}=A(t)x (t) , \ t\in\bbR, \  x(t) \in \bbR^d,
\end{equation}
where $A$ is a $d \times d$-matrix depending periodically on $t$.
Let $G(t)$ be the fundamental matrix-solution of (\ref{I1}). The Floquet theorem says that there exists a constant matrix $C$ and a periodic matrix $P(t)$
such that $G(t)=P(t)e^{Ct}$. This theorem allows to reduce the periodic system (\ref{I1}) to a system with constant coefficients.

We will consider the infinite dimensional case, where
$A$ is an operator in infinite dimensional Hilbert-spaces
and depends periodically on $t$. We will use the theory of this equation as developed and exposed in  \cite{Ku81}
and  \cite{KuchBook} (see in particular Sect. 5.1, the
comments of Ch.\,5 and the references there). This
allows us to avoid the details of the
Floquet-Bloch techniques, on which the theory is based
(for a recent presentation of the Floquet-Bloch techniques,
see  \cite{KuchRev}).
We also mention the related theory of elliptic problems in
periodic quasicylinders, see  \cite{na417} and \cite{NaPl},
although we stick to cylindrical domains here.

Our contribution to the infinite dimensional, periodic
problem will consist of a construction of a projection operator and a spectral splitting of the problem as described in the beginning. We also mention the
references \cite{KM1},  \cite{KM2} and  \cite{KMR},
the first two of which contain an analogous theory in the case
$A(t) $ is a perturbation of an operator $A_0$ independent of $t$. We will use the same formalism of analytic Fredholm operator in the way as it is presented in the appendix of
\cite{KM1}. In particular in the treatment of the infinite dimensional
part of the splitted system we use a technique developed \cite{KM1},
which  allows us to avoid a choice of function spaces for estimating
the remainder terms, since all of them can be treated from this
"pointwise estimate", see Sect.\,4 and 7 in \cite{KM1}.

The new aspects of the above approach are  motivated by the possibility
to generalize the existing center manifold reduction results of non-linear
Hilbert-space valued ODEs to the case of periodic equations. Indeed,
in the last section of the paper we will present such an application,
where we consider small solutions of a non-linear problem $D_tu-A(t)u=f(t,u)$
and present its reduction to a finite dimensional dynamical system, similarly
to the constant coefficient case considered by Mielke in \cite{M1}.

The structure of the paper is the following. In Sect.\,\ref{sec2} we
formulate the periodic infinite dimensional problem, introduce the function spaces and present the main assumptions
on the operator of the problem. In Sect.\,\ref{sec3} we remind some
basic definitions and properties of the eigenvalues, eigenvectors and
generalized eigenvectors of the operator pencils associated with our
periodic problem.  In Sect.\,\ref{sec4a} we collect known
results (\cite{KuchBook}, Ch.\,4 and 5) on the solvability and
asymptotics of solutions to periodic problems. These results are proved in
\cite{na417} in the case of elliptic boundary value problems with periodic coefficients in periodic cylinders. Our main results are contained in Sect.\,\ref{sec4} and \ref{sec5} . In Sect.\,\ref{sec4}
we introduce a certain  operator ${\mathcal P}$ which
is actually a projector which commutes both  with the periodic
operator as well
as the point evaluation operator $f \mapsto f(t)$ for all $t$,
and delivers a finite dimensional ODE with constant coefficients for ${\mathcal P}u$, where $u$ is a solution to the infinite dimensional system.  The main result of our paper is Theorem
\ref{T9ja} in Sect.\,\ref{sec5}, which gives a splitting of the
periodic operator valued equation into a finite dimensional system of ODEs  with constant coefficients and an infinite
dimensional part whose solutions have better properties in a certain sense.
In Section \ref{sec7} we present the main  application of Theorem \ref{T9ja},
namely an extension of the center manifold theorem for ODE in Hilbert space
obtained by Mielke in \cite{M1}, \cite{M}, to the periodic coefficient case. 

\section{Statement of the Problem}
\label{sec2}
The setting of our problem is similar to that of
\cite{KuchBook}, Section 5.1, and \cite{KM1}, Part III, except that
in the latter, the nature of the $t$-dependence of the operator
$A(t)$ is different. All the results of this section
are known, see in particular  Ch.\,1 of the first above mentioned
citation. Some proofs are reproduced for the convenience of the reader.

To proceed with
details, we let $X$ and $Y$  be Hilbert spaces over the complex scalar
field $\bbC$ such that $X$ is compactly and densely embedded in $Y$. We denote the norms in $X$ and $Y$  by $\Vert \cdot\Vert_X=\Vert \cdot ;X\Vert $ and $\Vert \cdot\Vert_Y=\Vert \cdot ;Y\Vert $,  respectively.  We identify $Y^*$ with $Y$ by using the inner product $(\cdot,\cdot)=(\cdot,\cdot)_Y$ and introduce for $h\in Y$ the norm
$$
\Vert h\Vert_{X^*}=\sup\{ |(g,h)|:g\in X, \Vert g\Vert_X=1\} .
$$
The completion of $Y$ with respect to this norm coincides with $X^*$, and the sesquilinear form $(g,h)$ can be extended for  $g\in X$ and $h\in X^*$ such
that the inequality $|(g,h)|\leq \Vert g\Vert_X\Vert h\Vert_{X^*}$
holds.  Clearly, $Y\subset X^*$.

Given $a,b \in \bbR$, $a < b$, we denote by ${\cX}(a,b)$   the space of functions
$u : (a,b) \mapsto X$ such that the  weak $t$-derivatives with values in $Y$
exist and are locally integrable (in the standard Bochner sense, see e.g. \cite{Hi}, Sect.\,3.7.) and such that the norm
\bea
& & \Vert u;{\mathcal X}(a,b)\Vert =\Big(\int\limits_a^b
\big(\Vert u(t);X\Vert ^2+\Vert D_t u(t);Y\Vert ^2 \big)dt\Big)^{1/2}  \label{1r}
\eea
is finite. Here and elsewhere $D_t=\partial/\partial t$. Also, the space
${\cY}(a,b)$ consists of locally integrable functions $u : (a,b) \mapsto Y$ with finite norm
\bea
\Vert f;{\mathcal Y}(a,b)\Vert =\Big(\int\limits_a^b\Vert f(t);Y\Vert ^2dt\Big)^{1/2}. \label{1s}
\eea
The space ${\mathcal Y}_{\rm loc}:=L^2_{\rm loc}(\bbR;Y)$ consists of
measurable functions defined on $\bbR$ with values in $Y$ with finite semi-norms
(\ref{1s}) for all $a<b$, and the space ${\mathcal X}_{\rm loc}:=
L^2_{\rm loc}(\bbR;X)$ is defined
analogously (cf. above); in particular for every $f \in {\mathcal X}_{\rm
loc}$, the semi-norms (\ref{1r}) are finite for all $a<b$ .

Given $\beta\in\bbR$, the space ${\mathcal X}_\beta$ consists of functions
$u\in L^2_{\rm loc}(\bbR;X)$ such that $
D_t u \in L^2_{\rm loc}(\bbR;Y)$ and
the norm
\begin{equation}\label{j6a}
\Vert u;{\mathcal X}_\beta\Vert =\Big(\int\limits_{\bbR}e^{2\beta t}(\Vert u(t);X\Vert ^2+\Vert D_t u(t);Y\Vert ^2)dt\Big)^{1/2}\ \mbox{is finite},
\end{equation}
and the space ${\mathcal Y}_\beta =L^2_\beta(\bbR;Y)$ consists of functions $f\in L^2_{\rm loc}(\bbR;Y)$ with  finite norm
\begin{equation}\label{j6b}
\Vert f;{\mathcal Y}_\beta\Vert =\Big(\int\limits_{\bbR}e^{2\beta t}\Vert f(t);Y\Vert ^2dt\Big)^{1/2}.
\end{equation}

In order to  deal with  periodic problems we follow \cite{N1}, \cite{NaPl} and also  introduce subspaces of
${\mathcal X}_{\rm loc}$ and ${\mathcal Y}_{\rm loc}$, which consist
of periodic functions in $t$ of period $1$ and which are  denoted by $\widehat{\mathcal X}$ and $\widehat{\mathcal Y}$, respectively. The norms in these spaces are
$$
\Vert u;\widehat{\mathcal X}\Vert =\Vert u;{\mathcal X}(0,1)\Vert ,\ \Vert f;\widehat{\mathcal Y}\Vert =\Vert f;{\mathcal Y}(0,1)\Vert .
$$

Let $A(t)$ be a bounded operator from $X$ into $Y$ depending continuously on
$t\in\bbR$ with respect to the operator norm. We assume  that $A(t)$ is
periodic with respect to $t$ with the period $1$.
For every $t$, we denote by $A (t)^* : Y  \to X^*$ the adjoint operator with respect
to the duality $(\cdot,\cdot)$, i.e.,
\bea
( A (t) \varphi, \psi )= ( \varphi, A(t)^* \psi )
\ \ \mbox \ \varphi \in X ,\psi \in Y . \label{1h}
\eea

We also define the  differential operators
\begin{equation}\label{1}
\cL= {\mathcal L}(t,D_t):=D_tu(t)+A(t)u(t) \  \ \mbox{and} \ \
{\mathcal L}^* (t,D_t):=-D_tu(t)+A(t)^* u(t) .
\end{equation}
In the following we will consider the problem
\begin{equation}\label{1k}
{\mathcal L}(t,D_t) u 
=f(t),
\end{equation}
where  $f\in L^2_{\rm loc}(\bbR;Y)$ is a given function and
$u\in {\mathcal X}_{\rm loc}$ is a function to be found. Our aim
is to introduce a reduction of this problem into a system
consisting of a  scalar valued, finite dimensional ODE-system and of another vector valued
ODE, which has better properties then the initial problem.
The first main assumptions on ${\mathcal L}$ is the following local estimate
(cf.  \cite{KM2}, Sect.\,2.2)
\bea
\Vert  u  ;  \cX (0,1)\Vert
& \leq &  C \Big(\Vert \mathcal{L}(t, D_t) u   ;  \cY (-1,2)   \Vert
\rowpl
\Vert  u   ;  \cY (-1,2)  \Vert\Big)
\ \ \ \ \   \mbox{for all } u \in     \cX (-1,2) .   \label{1a}
\eea

To formulate the second assumption
let us  consider the following operator depending on a complex parameter
$\lambda$,
\begin{equation}
\label{2}
{\mathcal A}(\lambda)={\mathcal L}(t,D_t)+\lambda:\widehat{\mathcal X}\rightarrow \widehat{\mathcal Y},\ \lambda\in\bbC .
\end{equation}
(Here, for the values $\lambda = i \xi$, the number
$\xi \in \bbR$ would correspond the Floquet parameter
or quasimomentum in the Floquet-Bloch transform, but we
do not need to exploit this machinery, as mentioned in the
introduction.)
Obviously, $\cA$  is a holomorphic operator
pencil with respect to the parameter $\lambda$. The second  main assumptions
on $\cL$ reads as  (cf.  \cite{N1}):
\bea
& &\mbox{\it  there exists $\lambda_0$ for which 
$\cA(\lambda_0) : \widehat{\mathcal X}\rightarrow \widehat{\mathcal Y}$ is an isomorphism.}   \label{2a}
\eea

{\bf Remark.}  We have in mind some applications to parabolic and elliptic  PDE-problems,
which have been transformed into first order ODE-systems with respect to
one of the variables  in a canonical way. The assumptions \ef{1a}, \ef{2a} are
natural for such cases. The assumptions would in general fail for  hyperbolic PDE-problems.

\BEL
\label{lem2.3}
If the assumptions \ef{1a} and \ef{2a} hold, then the families
\beas
& & \cA(\lambda) : \bbC \to \cL( \widehat \cX , \widehat \cY),  \ \ \
  \row
 \cA^*(\lambda):=\cA(\overline{\lambda})^*  : \bbC \to \cL( \widehat \cY  , \widehat \cX^* ),  \ \ \
\eeas
where $\cA^*(\lambda)=- D_t+A(t)^*+\overline{\lambda} $, are holomorphic Fredholm families.
Moreover, there holds
\begin{equation}\label{j3a}
\Big(\int\limits_0^1 \big( \Vert u (t);  X \Vert^2+\Vert (D_t+\lambda )u(t) ;  Y \Vert^2\big)dt\Big)^{1/2}\leq Ce^{|\Re\lambda|}\Big(\Vert {\mathcal A}(\lambda)u;\widehat \cY\Vert+\Vert u;\widehat \cY\Vert\Big)
\end{equation}
for $u\in \widehat \cX $.
\ENL

Here, $\cL( \widehat \cX , \widehat \cY)$ denotes the Banach space
of bounded linear operators from $ \widehat \cX$ into $\widehat \cY$.

\bigskip

Proof. Writing  estimate  \ef{1a} for the function $e^{\lambda t}u$, $u\in \widehat \cX $,  we get
\begin{eqnarray*}
&&\Big(\int\limits_0^1 e^{2t \Re\lambda }\big( \Vert u (t);  X \Vert^2+\Vert (D_t+\lambda )u (t);  Y \Vert^2\big)dt\Big)^{1/2}\\
&&\leq C\Big(\Big(\int\limits_{-1}^2e^{2t \Re\lambda}\Vert  \cL(t, D_t+\lambda) u(t);Y\Vert ^2dt\Big)^{1/2}+
\int\limits_{-1}^2 e^{2t \Re\lambda }\Vert  u(t);Y\Vert ^2dt\Big)^{1/2}\Big),
\end{eqnarray*}
which implies  estimate (\ref{j3a}).
Since the inclusion $X\subset Y$ is compact, we can use the argument
in \cite{LM}, p. 20 or Theorem 2.1. to see that the embedding $\widehat \cX
\subset\widehat \cY$ is also compact. Hence,   estimate (\ref{j3a}) implies that
the kernel of ${\mathcal A}(\lambda)$ is finite dimensional and the image is
closed for all $\lambda$. This together with assumption (\ref{2a}) gives that
the operator pencil is Fredholm with the index $0$ for all $\lambda$ (see
\cite{KM1}, Section A.8).

The definition of the adjoint holomorphic family $\cA^* (\lambda)$ is as in
\cite{KM1}, Section A.9, and its Fredholm property follows from that of
the family $\cA(\lambda)$, as explained in the citation; see also the next lemma.  \ \
 $\Box$

\bigskip

Let us provide a description of the dual space $\widehat \cX^*$. To this end we use
the inner product to  identify $\widehat \cY^* =\widehat \cY$, and we also
denote by $L^2_{\rm per}(\bbR;X^*)$  the subspace of $L^2_{\rm loc}
(\bbR;X^*) $ consisting of periodic functions $f : \bbR \to X^*$, endowed
with the norm
\beas
\Vert f ;L^2_{\rm per}(\bbR;X^*) \Vert
= \Big( \int\limits_0^1 \Vert f(t) ; X^* \Vert^2 dt \Big)^{1/2}.
\eeas
We skip the standard proof of the  following lemma.

\BEL
\label{lem2.1}
Under the dual pairing
\beas
& & ( u , v )_{\widehat \cY} :=
\int\limits_0^1  ( u(t), v(t) )  \, dt \ , \ \
\eeas
the dual space $\widehat \cX^* $ of $\widehat \cX$
consists of periodic functions $w$ represented as
\bea
w=w_0+D_tw_1 ,  \label{j4p}
\eea
where $w_0\in L^2_{\rm per}(\bbR;X^*)$ and $w_1\in \widehat{\mathcal Y}$,
and it is endowed with the norm
\beas
& & \Vert w;\widehat{\mathcal X}^* \Vert =\inf \big(\Vert w_0;L^2_{\rm per}(\bbR;X^*) \Vert+\Vert D_t w_1; \widehat{\mathcal Y}\Vert \big) \ , \ \
\eeas
where the infimum is taken over all representations \ef{j4p}.
The adjoint $\cA^*(\lambda)$ of the operator $\cA(\lambda)$ satisfies
\beas
( \cA(\lambda)  \varphi, \psi )_{\widehat \cY} = ( \varphi,
\cA^* (\lambda) \psi )_{\widehat \cY}
\ \ \mbox \ \varphi \in \widehat \cX ,\psi \in \widehat \cY . 
\eeas
\ENL

We end this section with one more lemma.  One can easily verify that
\begin{equation}\label{j3c}
e^{2\pi it} \cA(\lambda)u= \cA( \lambda-2\pi i)(e^{2\pi it}u)
\end{equation}
for $u\in \widehat \cX $.

\BEL
\label{Lj3} Assume that the operator ${\mathcal A}(\mu)$ is an isomorphism for $\mu=\beta+i\xi$ with a fixed $\beta \in \bbR$
and for all $\xi\in [0,2\pi)$. Then, for all $u\in \widehat \cX $
and $\lambda=\beta+i\xi$  with  $\xi\in\bbR$,
\begin{equation}\label{j3d}
\Big(\int\limits_0^1  \big(\Vert u (t);  X \Vert^2+\Vert (D_t+\lambda )u(t) ;
Y \Vert^2\big) dt\Big)^{1/2}\leq C\Vert {\mathcal A}(\lambda)u;
\widehat \cY\Vert,
\end{equation}
where $C$ may depend on $\beta$ but it is independent of $\xi$.
\ENL
Proof. By (\ref{j3c}) the optimal constant $c$ in the inequality
\begin{equation}\label{j3e}
\Vert \cA(\lambda)u;\widehat \cY\Vert\geq c\Vert u;\widehat \cY\Vert
\end{equation}
is the same for $\lambda$ and $\lambda-2\pi k i$ for all $k
= \pm 1, \pm 2 , \ldots$. This together with (\ref{j3c}) and the assumption of the lemma implies existence of a constant $c_0$ such that (\ref{j3e}) is true for all $\lambda=\beta+i\xi$  with $\xi\in\bbR$. Using (\ref{j3e}) we derive (\ref{j3d}) from (\ref{j3a}).
 \ \ $\Box$

\bigskip

\section{Eigenvectors, generalized eigenvectors, Jordan chains}
\label{sec3}
We recall some basic facts concerning the spectrum of the operator
pencil  $\cA(\lambda)$, \ef{2}. For the sake of the presentation, we
recall some proofs and otherwise refer to the Appendix \cite{KM1}
for a quick introduction to the topic. As in standard spectral theory of linear operators,
the spectrum is  the set of those $\lambda
\in \bbC$ such that $\cA(\lambda) : \widehat \cX \to \widehat \cY$ is not
invertible;  $\lambda$ is an eigenvalue, if the kernel
of $\cA(\lambda)$ is not $\{ 0 \}$.

Since
$\cA(\lambda) : \widehat{\mathcal X}\rightarrow \widehat{\mathcal Y} $
is a holomorphic Fredholm family and due to the assumption (\ref{2a}), the spectrum of $\cA(\lambda)$ consists of isolated eigenvalues
of finite algebraic multiplicity, see Proposition A.8.4 of \cite{KM1}.
From the relation (\ref{j3c})
it follows that if $\lambda$ is an eigenvalue then the same is true for
$\lambda+2\pi i$ and their multiplicities coincide. In the following we denote
for all $\beta \in \bbR $
\begin{equation}\label{4}
\delta_\beta=\{\lambda\in\bbC\,:\, \Re\lambda = \beta,\,\Im\lambda\in [0,2\pi )\} ,
\end{equation}
and we  choose real numbers
\bea
\beta_1 <  \beta_2 \label{4ab}
\eea
such that there are no eigenvalues of $\cA (\lambda) $ on the the intervals $\delta_{\beta_1}$ and $\delta_{\beta_2}$.
We denote  eigenvalues of ${\mathcal A}(\lambda)$ in the set
\bea
\{ \lambda=\beta+i\xi\,;\,\beta_1<\beta<\beta_2,\;\xi\in[0,2\pi)\}
\label{4abx}
\eea
 by $\lambda_k$, $k=1,\ldots,N$, and let $J_k$ and $m_{k,1},\ldots,m_{k,J_k}$ be the geometric and partial multiplicities  of $\lambda_k$.
Assume that for every $k=1, \ldots, N$,
\bea
\varphi^k_{j,m},\;m=0,\ldots,m_{k,j}-1,\;j=1,\ldots,J_k,   \label{4b}
\eea
is a canonical system of Jordan chains of the linear pencil $\cA( \lambda) $ corresponding to $\lambda_k$ (see \cite{KM1}, Definition A.4.3,
Propositions A.4.4, A.4.5.). 
The functions
\bea
\varphi^k_{j,0 }, \  j=1,\ldots,J_k,  \label{4b3}
\eea
form a linearly independent sequence of eigenvectors corresponding to the eigenvalue $\lambda_k$, while the functions
\ef{4b} with $m \geq 1$ are associated vectors satisfying
\bea
\cA(\lambda_k)  \varphi_{j,0 }^k=0,\ \;\cA(\lambda_k)  \varphi_{j,m }^k =   -\varphi_{j,m- 1 }^k   , \ \ \
m= 1, \ldots , m_{k,j} -1 .
\label{4d}
\eea
In the same way, the eigenfunctions and generalized eigenfunctions of the
adjoint operator are the solutions of the equations
\bea
\cA^* ( {\lambda_k})\psi_{j,0 }^k=0,\ \;\cA^* ({\lambda_k})\psi_{j,m }^k =  -\psi_{j,m -1 }^k , \ \ \
m= 1, \ldots , m_{k,j} -1 .
\label{4dd}
\eea
It will be important to specify the choice of the functions  \ef{4b3}
and \ef{4dd} such that certain orthogonality relations are satisfied.
Notice that we consider a finite set of  eigenvalues, which is fixed
by the choice of the numbers $\beta_1$, $\beta_2$ above. The following assertion is known and its proof can be found in
  Remark A.10.3 in \cite{KM1} (see formula (A.60) there).

\BEL
\label{lem3.0} If the Jordan chains (\ref{4b}) are fixed, then there exist
uniquely defined  Jordan chains of the adjoint pencil $\cA^*(\lambda)$ corresponding to the eigenvalue $\overline{\lambda_k}$
\bea
\psi^k_{j,m},\;m=0,\ldots,m_{k,j}-1,\;j=1,\ldots,J_k,   \label{4bb}
\eea
such that in addition to all equations \ef{4d} and \ef{4dd} also the
following hold true:
\bea
( \varphi_{j,m_{k,j}-1}^k ,  \psi_{J,m}^k   )_{
\widehat \cY} = \delta_j^J\delta_0^m,\,\,m=0,\ldots,m_{k,j}-1\,.   \label{4ddd}
\eea
\ENL

\bigskip

From now on we assume that the eigenfunctions and generalized
eigenfunctions satisfy \ef{4ddd}. The last relation implies some more orthogonality relations.

\BEL
\label{lem2.5} Let the Jordan chains \ef{4b} and \ef{4bb} be the same as in Lemma \ref{lem3.0}. Then the following biorthogonality
relations hold:
\bea
( \varphi_{j, m}^k , \psi_{J,m_{K,J}-1 -M}^K )_{\widehat \cY}
= \delta_k^K \delta_j^J \delta_m^M
\label{4de}
\eea
for all $k,K,j,J,m,M$.
\ENL

Proof. Let first $K=k$. Then (\ref{4de}) for $m=m_{k,j}-1$ follows from (\ref{4ddd}).

Next we observe that for $ m = 1 , \ldots ,
m_{k,j}-1 $ and $M = 1 , \ldots , m_{k,J} -1$, $J=1,\ldots,J_k, $ the relations  \ef{4d}
and  \ef{4dd} yield
\bea
& & ( \varphi_{j, m}^k , \psi_{J,m_{k,J}-1-M }^k )_{
\widehat \cY} =
 - ( \varphi_{j, m}^k , \cA^*({\lambda_k})  \psi_{J,m_{k,J}-M }^k )_{
\widehat \cY}
\roweq
- ( \cA(\lambda_k)  \varphi_{j, m}^k , \psi_{J,m_{k,J}-M }^k )_{
\widehat \cY} =
( \varphi_{j, m -1 }^k , \psi_{J,m_{k,J}-M }^k )_{
\widehat \cY}. \label{4dj}
\eea
Applying this relation with $m=m_{k,j}-1$ and $M=1,\ldots,m_{k,J}-1$ and using that (\ref{4de}) is proved for $m=m_{k,j}-1$, we arrived at (\ref{4de}) for $m=m_{k,j}-2$ and $M=0,\ldots, m_{K,J}-2$. Since the relations
$$
( \varphi_{j, m}^k , \psi_{J,0}^k )_{\widehat \cY}=0,\;\; m=0,\ldots,m_{k,j}-2,
$$
follow from the solvability of (\ref{4d}), we arrive at (\ref{4de}) for $m=m_{k,j}-2$ and all $M$.
Repeating this argument we prove (\ref{4de}) for all $m$ and $M$.

We finally show that if $k \not= K$, then the orthogonality in \ef{4de}
automatically holds. For the two eigenfunctions we get the orthogonality
$(\varphi_{j,0}^k , \psi_{J,0}^K )_{\widehat \cY} = 0$ for all
$j  = 1, \ldots , J_k$, $J = 1, \ldots , J_K$ by the
simple classical argument, since the eigenvalues $\lambda_k$ and
$\lambda_K$ are different. Then, we have, for all $M= 0 , \ldots , m_{K,J}-
2$, all $j, J$,
\bea
& & (\varphi_{j,0}^k , \psi_{J,M}^K )_{\widehat \cY}
=  - ( \varphi_{j,0}^k , \cA^* ({\lambda_K}) \psi_{J,M+1}^K )_{\widehat \cY}
=  - ( \cA (\lambda_K) \varphi_{j,0}^k , \psi_{J,M+1}^K )_{\widehat \cY}
\roweq
- (\cA (\lambda_k)  \varphi_{j,0}^k ,  \psi_{J,M+1}^K )_{\widehat \cY}
+ (\lambda_k- \lambda_K) ( \varphi_{j,0}^k ,  \psi_{J,M+1}^K )_{\widehat \cY}
\roweq
(\lambda_k- \lambda_K) ( \varphi_{j,0}^k ,  \psi_{J,M+1}^K )_{\widehat \cY} , \label{4rr}
\eea
where the coefficient $\lambda_k- \lambda_K$ is non-zero,
so that the orthogonality $( \varphi_{j,0}^k ,  \psi_{J,M+1}^K )_{\widehat
\cY} = 0$ for all $M = 0 , \ldots , m_{K,J}- 1$ and $j,J$ follows by
induction. In the same way one obtains $( \varphi_{j,m}^k ,  \psi_{J,0}^K
)_{\widehat  \cY} = 0$ for all $m = 0 , \ldots , m_{k,j}- 1$ and $j,J$ .

Then, one proves the following formulas in the same way as \ef{4rr}
\begin{equation*}
(\varphi_{j,m}^k , \psi_{J,M}^K )_{\widehat \cY}
=  \left\{
\begin{array}{rr}
(\varphi_{j,m+1}^k , \psi_{J,M-1 }^K )_{\widehat \cY}
+ (\lambda_k - \lambda_K) ( \varphi_{j,m+1}^k ,  \psi_{J,M}^K )_{\widehat \cY} & \\
& \\
(\varphi_{j,m-1}^k , \psi_{J,M+1 }^K )_{\widehat \cY}
+ (\lambda_k- \lambda_K) ( \varphi_{j,m}^k ,  \psi_{J,M+1}^K )_{\widehat \cY} \ .&
\end{array}
\right.
\end{equation*}
One can then proceed by induction 
to get the orthogonality for all indices. \ \ $\Box$

\bigskip

Let us still introduce some more notation with the help of the above
introduced Jordan chains: we  define
\begin{equation}\label{j6d}
\Phi_{j,m}^k (t)= e^{\lambda_k t}\sum_{n=0}^{m}\frac{t^n}{n!}\varphi^k_{j,m-n}
= e^{\lambda_k t}\sum_{n=0}^{m}\frac{t^{m-n}}{(m-n)!}\varphi^k_{j,n},
\end{equation}
for all $k=1,\ldots,N, \ j=1,\ldots,J_k, \ m=0,\ldots,m_{k,j}-1.$
It is known and one can verify it directly that these functions are solutions to
the homogeneous equation (\ref{1k}). The binomial formula implies the
following relation which will be needed later:
\bea
& & e^{\lambda_kt}\sum_{n=0}^m\frac{(t-\tau)^n}{n!}\varphi_{j,m- n}^k(t)
=
e^{\lambda_kt}\sum_{n=0}^m \sum_{\nu = 0}^n
\frac{t^{n - \nu} (-\tau)^\nu }{(n-\nu)!\nu!}\varphi_{j,m- n}^k(t)
\roweq
e^{\lambda_kt}\sum_{\nu=0}^m   \frac{(-\tau)^\nu}{\nu!} \sum_{n = 0}^{m-\nu}
\frac{t^{m - \nu -n }  }{(m -\nu -  n)!}\varphi_{j,n}^k(t) = 
\sum_{\nu=0}^m\frac{(-\tau)^{\nu}}{\nu!}\Phi_{j,\nu}^k(t).
\label{rela}
\eea

\section{Some results on solvability and asymptotics for problem (\protect\ref{1k})}
\label{sec4a}

We will need some more solvability and asymptotical results for
problem \eqref{1k}, which are proved for general boundary value problems  with periodic coefficients in a cylinder in \cite{KuchBook},
Sect. 4.2, 5.1, 5.4 (for the case of
a periodic quasi-cylinder, see also \cite{na417}).


\begin{theorem}\label{Tj6a} The mapping
\beas
{\mathcal L}(t,D_t):{\mathcal X}_{-\beta}\rightarrow {\mathcal Y}_{-\beta}
\eeas
is isomorphic if 
the semi-interval $\delta_{ \beta}$ does not contain eigenvalues of the
operator pencil ${\mathcal A}(\lambda)$.
\end{theorem}

\BET\label{Tj6b}
Let $\beta_1<\beta_2$ be real numbers such that the semi-intervals $\delta_{\beta_1}$ and $\delta_{\beta_2}$ do not contain eigenvalues of the operator pencil ${\mathcal A}(\lambda)$ and let $f\in {\mathcal Y}_{-\beta_1}\bigcap {\mathcal Y}_{-\beta_2}$. Denote by $u_1$ and $u_2$
solutions to the problem {\rm (\ref{1k})} from the spaces ${\mathcal X}_{-\beta_1}$ and ${\mathcal X}_{-\beta_2}$ respectively (which exist according to
Theorem \ref{Tj6a}). Then
\begin{equation}\label{j6e}
u_2-u_1=\sum_{k=1}^{N}\sum_{j=1}^{J_k}\sum_{m=0}^{m_{kj}-1}c_{j,m}^k
\Phi_{j,m}^k (t),
\end{equation}
where the functions $\Phi_{j,m}^k $, $k = 1, \ldots ,N$,  are all the functions (\ref{j6d}) such that  the eigenvalues $\lambda_k$  belong to  the set
\bea
 Q(\beta_1, \beta_2) := \{\lambda = \beta + i \xi \, : \, \beta_1
 <  \beta < \beta_2 , \ \xi \in [0,2\pi) \} ,    \label{59o}
\eea
and   $c_{j,m}^k$ are constants.
\ENT

\bigskip

A straightforward consequence of the last theorem is the following uniqueness result.

\begin{corollary}\label{Cor1a} Let $\beta_1$ and $\beta_2$ be the same as in Theorem \ref{Tj6b}. If $u\in {\mathcal X}_{\rm loc}$ is a solution of (\ref{1k}) with $f=0$ and
\begin{equation}
\label{j7b}
\Vert u;{\mathcal X}(t,t+1)\Vert \leq Ce^{\beta_1 t}\;\,\;\mbox{for $t\geq 0$ and }\ \Vert u;{\mathcal X}(t,t+1)\Vert \leq Ce^{\beta_2 t}\;\,\;\mbox{for $t\leq 0$}
\end{equation}
for some positive constant $C$, then
\begin{equation}\label{j7a}
u=\sum_{k=1}^{N}\sum_{j=1}^{J_k}\sum_{m=0}^{m_{k,j}-1}c_{j,m}^k
\Phi_{j,m}^k (t),
\end{equation}
where $c_{j,m}^k $ are constants and  $\Phi_{j,m}^k $ are all
functions (\ref{j6d}) such that the eigenvalues $\lambda_k$ belong to the set
\eqref{59o}.
\end{corollary}

Proof. Let $\beta'_1 > \beta_1$ and $\beta_2' < \beta_2$ be such that the intervals $[\beta_1,\beta_1']$ and $[\beta_2',\beta_2]$ do not contain the eigenvalues of the operator pencil
${\mathcal A}(\lambda)$. Let also $\eta=\eta(t)$ be a smooth function of one variable such that $\eta(t)=1$ for $t>1$ and $\eta(t)=0$ for $t>0$. Consider the problem (\ref{1k}) with $f= (D_t\eta)u$. This problem has two solutions,  $u_1=\eta u$ and $u_2=(\eta-1)u$. Since
\begin{eqnarray*}
&& \int\limits_{0}^\infty e^{-\beta_1't}\Vert u_1(t);X\Vert ^2dt
\leq C\int\limits_{0}^\infty e^{-\beta_1' t }\Vert u;{\mathcal X}(t,t+1)\Vert ^2dt
\leq 
C\int\limits_{0}^\infty e^{t(\beta_1-\beta_1')}dt<\infty,
\end{eqnarray*}
$u_1$ belongs to $\cX_{-\beta_1'}$. Similarly, $u_2$ belongs  to $\cX_{-\beta_2'}$. By Theorem \ref{Tj6b},
$u=u_2-u_1$ is equal to the right-hand side of (\ref{j6e}) and we arrive at
(\ref{j7a}). \ \ $\Box$

\section{Pointwise projector.} \label{sec4}

Let us introduce the operator 
\begin{equation}
\label{L1}
{\mathcal P}u(t)=\sum_{k=1}^N\sum_{j=1}^{J_k}\sum_{m=0}^{m_{k,j}-1}
u^k_{j,m  }(t)\varphi^k_{j, m_{k,j}-1-m}(t),
\end{equation}
where
\begin{equation}\label{KK1}
u_{j,m  }^k (t) = \big( u(t) , \psi_{j,m  }^k(t) \big)
\end{equation}
and $t \in\bbR$. If the inner product $(\cdot,\cdot)$ were replaced by
$(\cdot,\cdot)_{\widehat{\mathcal Y}}$ then by (\ref{4de}) the above
operator would be the finite dimensional spectral projector for the
operator $-{\mathcal L}(t,D_t)$ corresponding to the spectrum  in the set
(\ref{4abx}). But here we have only the  inner product $(\cdot,\cdot)$ on
the cross-section and the coefficients $u^k_{j,m  }$ are functions of $t$.
By definition, the operator $\cP$ commutes
with the point evaluation operator $ \cY_{\rm loc} \ni f \mapsto f(t)$
for all $t \in \bbR$. The notation ${\mathcal P}(t)$ means the evaluation
of \eqref{L1} at $t$.

The operator ${\mathcal P}$ is well defined on the spaces
${\mathcal Y}_{\rm loc}$ and ${\mathcal X}_{\rm
loc}$. Since there are only finitely terms in the
sums \eqref{L1} and
$$
|u^k_{j,m}(t)|\leq C\Vert u(t);Y\Vert \ \mbox{and}\ |\partial_tu^k_{j,m}(t)|
\leq C\Vert \partial_tu(t);Y\Vert ,
$$
the operator ${\mathcal P}$ is bounded as an operator in ${\mathcal X}_\beta$ and ${\mathcal Y}_\beta$.

In the following theorem we prove that the operator ${\mathcal P}(t)$ is
indeed a projector and derive some important properties of ${\mathcal P}$.

\begin{theorem}
\label{Tmain} The operator ${\mathcal P}$ has the following properties:

\noindent (i) ${\mathcal P}(t)$ is a projector for every $t\in \bbR$;

\noindent (ii) ${\mathcal L}{\mathcal P}u={\mathcal P}{\mathcal L}u$ for all $u\in {\mathcal X}_{\rm loc}$ ;

\noindent (iii) there holds
\begin{equation}\label{10a}
{\mathcal L}{\mathcal P}u=\sum_{k=1}^N\sum_{j=1}^{J_k}\sum_{m=0}^{m_{k,j}-1}\Big((D_t-\lambda_k)u^k_{j,m}(t)-u^k_{j,m-1}(t)\Big)\varphi^k_{j,m_{k,j}-1-m}(t),
\end{equation}
where it is assumed that $u^k_{j,-1}=0$.

\end{theorem}

\bigskip
The proof of is based on the following, perhaps unexpected fact.

\BEL
\label{lem2.6}
We have for all $k,K =1, \ldots , N$, $j=1, \ldots, J_k$,
$J=1, \ldots, J_K$,  $m=0, \ldots, m_{k,j}-1$,
$M=0, \ldots, m_{K,J}-1$,
\begin{equation}\label{L2}
\big( \varphi^k_{j,m_{k,j}-1-m}(t),\psi^K_{J,M}(t) \big)=\delta^k_K\delta^j_J\delta^m_M \ \ \ \mbox{for all $t\in\bbR$}.
\end{equation}
\ENL

We will give a proof of this lemma at the end of this section and now
complete the proof of Theorem \ref{Tmain}. Indeed, the claim $(i)$ follows
from   (\ref{L2}). As for $(iii)$, using (\ref{4d}) we have
\begin{eqnarray*}
&&(D_t+A(t))U(t)
\roweq
\sum_{k=1}^N\sum_{j=1}^{J_k}\sum_{m=0}^{m_{k,j}-1}
\Big((D_tu^k_{j,m}(t))\varphi^k_{j,m_{k,j}-1-m}(t)  \\
& & \hskip2.5cm + \,
u^k_{j,m}(t)(D_t+A(t))\varphi^k_{j,m_{k,j}-1-m}(t)\Big)
\roweq
\sum_{k=1}^N\sum_{j=1}^{J_k}\sum_{m=0}^{m_{k,j}-1}
\Big((D_t-\lambda_k)u^k_{j,m}(t)\varphi^k_{j,m_{k,j}-1-m}(t)
 \\
& & \hskip2.5cm - \,
u^k_{j,m}(t)\varphi^k_{j,m_{k,j}-m-2}(t)\Big),
\end{eqnarray*}
which implies $(iii)$.
Finally, we have
\beas
& & {\mathcal P}{\mathcal L}u=\sum v^k_{j,m}(t)\varphi^k_{j,m_{k,j}-1-m}(t) \ , \ \ \mbox{where}
\row
v^k_{j,m}(t)=\big( (D_t+A(t))u(t) , \psi_{j,m  }^k(t) \big)=D_tu^k_{j,m}(t)+\big( u(t) , (-D_t+A(t)^*)\psi_{j,m  }^k(t) \big).
\eeas
Using (\ref{L2}) we obtain
\beas
v^k_{j,m}(t)&=&
D_tu^k_{j,m}(t)-\lambda_k\big( u(t) , \psi_{j,m  }^k(t) \big)-\big( u(t) ,
\psi_{j,m-1  }^k(t) \big)
\roweq
(D_t-\lambda_k)u^k_{j,m}(t)-u^k_{j,m-1}(t).
\eeas
which together with $(iii)$ gives $(ii)$.
 \ \ $\Box$

\bigskip

Proof of Lemma \ref{lem2.6}.  Introduce
\bea
I(k,j,m;K,J,M)(t)= \big( \varphi^k_{j,m}(t),\psi^K_{J,M}(t) \big).
\label{L2y}
\eea
 We have
for all $m,M=0, \ldots ,  m_{k,j}-1$ 
\begin{eqnarray*}
&& D_tI(k,j,m;K,J,M)(t)=
(D_t\varphi^k_{j,m}(t),\psi^K_{J,M}(t))+(\varphi^k_{j,m}(t),D_t\psi^K_{J,M}(t)).
\end{eqnarray*}
Moreover, by \ef{4d}, \ef{4dd},
$$
(D_t+A(t)+\lambda_k )\varphi^k_{j,m}+\varphi^k_{j,m-1}=0,\ (-D_t+A(t)^*+\overline{\lambda_K})\psi^K_{J,M}+\psi^K_{J,M-1}=0,
$$
where we must agree that the functions with negative second lower index are
zero. Thus,
\begin{eqnarray*}
D_tI(k,j,m;K,J,M)(t)=
 &&-\big((A(t)+\lambda_k)\varphi^k_{j,m}(t)+
\varphi^k_{j,m -1}(t),\psi^K_{J,M}(t)\big)\\
&&+\big(\varphi^k_{j,m}(t),(A(t)^*
+\overline{\lambda_K})\psi^K_{J,M}(t)+\psi^K_{J,M-1}(t)\big).
\end{eqnarray*}
After cancellation, we get
\bea
\label{L3a}
& & D_tI(k,j,m;K,J,M)(t)
\roweq
-\big(\varphi^k_{j,m -1 }(t),\psi^K_{J,M}(t) \big)
+\big(\varphi^k_{j,m}(t),\psi^K_{J,M-1}(t)\big)
\roweq
- I(k,j,m - 1;K,J,M)(t) + I(k,j,m;K,J,M-1)(t).
\eea

$1^\circ.$ We first prove the case
$K=k $ and $J = j \in \{1, \ldots, J_k \}$. 

$(i)$ Let first  $m = 0 $ and  $M=0$.
Then,  the right-hand side of \ef{L3a} is zero by the
convention made above, and therefore
$D_tI(k,j,0;k,j,0)(t) = 0$ and thus  $I(k,j,0;k,j,0)(t)$ does not depend on $t$. By Lemma \ref{lem2.5} we get for all $t$
\bea
& & I(k,j,0;k,j,0)(t) = \int\limits_0^1 I(k,j,0,;k,j,0)(\tau) d \tau
= 
( \varphi_{j, 0 
}^k , \psi_{j,0 }^k )_{\widehat \cY}
= 0 . 
\label{L3b} 
\eea

$(ii)$ We next consider the case $m = 0 $,
$M= 1, \ldots , m_{k,j} -2$. We use  induction with respect to $M$: assume that
$I(k,j,0 ;k,j,M) (t) = 0$ for some $M < m_{k,j} -2$ and all $t \in
[0,1]$. By \ef{L3a} we get
\bea
& &  D_t I(k,j, 0;k,j,M +1) (t)
= 
I(k,j,0;k,j,M) (t) = 0  \label{L3c}
\eea
hence, $I(k,j,0 ;k,j,M +1)(t)$ is constant with respect to $t$.
Integrating this constant as in \ef{L3b} and using Lemma \ref{lem2.5}
yield for all $t$
\beas
& & I(k,j,0;k,j,M) (t)
= 
(\varphi_{j, 0}^k ,\psi_{j,M}^k )_{
\widehat{\cY}}  = 0 \ \ \forall \,
M= 1, \ldots , m_{k,j} -2. 
\eeas
(Note that by Lemma \ref{lem2.5} the inner product is not zero
for $M=m_{k,j}-1$.)

In the same way, using inductively
\beas
D_t I(k,j,m +1 ;k,j, 0 ) (t) =
I(k,j,m  ;k,j,0) (t) = 0 
\eeas
for $m= 0, \ldots , m_{k,j}-2 $ instead of \ef{L3c} we prove that
\beas
I(k,j,m;k,j, 0 ) = 0 \ \ \forall \,
m = 1, \ldots , m_{k,j}-2 . 
\eeas

$(iii)$ We next consider the case $m + M \leq m_{k,j} - 2 $ by using
a double induction: assume that for some index $M\geq 0$ with
$ M \leq m_{k,j} - 3 $ 
the equality $I(k,j,m;k,j,M) = 0$ has been proven for all
$m = 0, \ldots ,  m_{k,j} -  M - 3 $.
Then,  \ef{L3a} implies, for $m = 1,  \ldots,  m_{k,j} -  M - 2,$
\beas
& &  D_t  I(k,j,m ;k,j,M+1) \roweq
- I(k,j,m-1;k,j,M+1) + I(k,j,m ;k,j,M) .
\eeas
We can thus proceed by induction with respect to $m$  (using $t$-integration
and Lemma \ref{lem2.5} as above) to get
\beas
I(k,j,m;k,j,M+1) = 0  \ \ \forall  \,
m = 0, \ldots , m_{k,j}- M-2.
\eeas
Induction with respect to $M $ yields \ef{L2}
for all $m + M \leq m_{k,j} - 2$. In both induction procedures we
use $(ii)$ for $m=0$ and $M=0$ to start with.

$(iv)$ Consider   $m + M = m_{k,j} - 1$. Formula \ef{L3a} and what we have
proven until now again imply that for every  $m= 0,1 \ldots, m_{k,j} -1$ the
expression
$$
I(k,j,m;k,j, m_{k,j} - m -1 )
$$ is a constant, which by $t$-
integration and Lemma \ref{lem2.5} is equal to 1.

$(v)$
However, this and \ef{L3a}
again imply that for every $m= 1,2,\ldots, m_{k,j} -1$ we have
\beas
I(k,j,m;k,j,  m_{k,j} - m ) = 0   .
\eeas

$(vi)$
From here on we can continue in the same way as in $(iii)$ to get the
result for $M + m \geq  m_{k,j}$.




\bigskip

$2^\circ.$ The proof in the case $K=k$ but $J\not= j$ is simpler
than $1^\circ$. The case $(i)$ is the same. The argument of the case $(ii)$
yields  $I(k,j,m;k,J,M) =0$  for the pairs $(m,M)$ with
$m=0, M=1, \ldots , m_{k,J}-1$ and $m=1, \ldots , m_{k,j}-1$, $M=0$,
by Lemma \ref{lem2.5}. Then, the procedure of $(iii)$ yields
$I(k,j,m;k,J,M) =0$ for all remaining pairs $(m,M)$, since  now we
do not have the obstruction of the case $(iv)$ (the inner products
$(\varphi_{j,m}^k, \psi_{J,M}^k )_{\widehat \cY} $ equal 0 instead of
1, by Lemma \ref{lem2.5}, for all indices in question).

$3^\circ.$ For the proof in the case $K\not=k$ we need to introduce
instead of \ef{L2y},
$$
I(k,j,m;K,J,M)(t)=e^{(\lambda_k-\lambda_K)t}
\big( \varphi^k_{j,m}(t),\psi^K_{J,M}(t) \big),
$$
because $\lambda_k \not= \lambda_K$. By a similar calculation as
around \eqref{L3a} we get for all $j,J,m,M$
\bea
&& e^{(-\lambda_k+\lambda_K)t} D_tI(k,j,m;K,J,M)(t)
\roweq
(\lambda_k-\lambda_K)\big( \varphi^k_{j,m}(t),\psi^K_{J,M}(t) \big)
 \rowpl
(D_t\varphi^k_{j,m}(t),\psi^K_{J,M}(t))+(\varphi^k_{j,m}(t),
D_t\psi^K_{J,M}(t))
\roweq
- \big( \varphi^k_{j,m-1}(t),\psi^K_{J,M}(t) \big)
+ \big( \varphi^k_{j,m}(t),\psi^K_{J,M-1}(t) \big) .
\label{L3at}
\eea
The following argument shows that we can use \ef{L3at} instead of
\ef{L3a} and repeat the proof of the case $2^\circ$ (i.e. the steps
$(i)-(iii)$ in $1^\circ$) for all $j,J,m,M$:
assume that the right-hand side of \ef{L3at} equals 0 for all $t$.
Since $ e^{(-\lambda_k+\lambda_K)t} \not=0$, this implies
that  $ D_tI(k,j,m;K,J,M)(t) = 0$  for all $t$, hence,
$ I(k,j,m;K,J,M)(t) = B$ for some constant $B$, for all $t$. Thus, by Lemma \ref{lem2.5},
\beas
& & 0 = \big( \varphi^k_{j,m},\psi^K_{J,M} \big)_{\widehat \cY}
= \int\limits_0^1 e^{(\lambda_K-\lambda_k)t}
e^{(\lambda_k-\lambda_K)t}
\big( \varphi^k_{j,m}(t),\psi^K_{J,M}(t) \big) dt
= B \int\limits_0^1 e^{(\lambda_K-\lambda_k)t} dt .
\eeas
Here we have $\int_0^1  e^{(\lambda_K-\lambda_k)t} dt \not =0$,
since $\lambda_K-\lambda_k$ cannot equal a multiple of $i 2 \pi$,
see \eqref{4abx}.
Hence, the constant $B$ must be zero, and thus also
$\big( \varphi^k_{j,m}(t),\psi^K_{J,M}(t) \big) = 0$  for all $t$.
  \ \ $\Box$

\bigskip

\section{Spectral splitting.}\label{sec5}

We assume that $\beta_1 < \beta_2$ are the same real numbers as in
Theorem \ref{Tj6b} so that in particular the semi-intervals
$\delta_{\beta_1}$, $\delta_{\beta_2}$ do not contain eigenvalues
of the pencil $\cA(\lambda)$ and its eigenvalues in the set
$Q(\beta_1, \beta_2)$ of \ef{59o} are $\lambda_1, \ldots, \lambda_N$.
We introduce the function
\beas
\mu (t)=e^{-\beta_1t}\ \mbox{for $t\geq 0$ and}\ \mu (t)=e^{-\beta_2t}\ \mbox{for $t\leq 0$.}
\eeas

Using the projector $\cP$ of \ef{L1} we represent a solution of (\ref{1k}) as
$$
u=U+V,\ \  \mbox{where} \ U={\mathcal P}u, \ V={\mathcal Q}u
:= ( I-{\mathcal P}) u.
$$
Due to the commutation relation Theorem \ref{Tmain}(ii) we have the following system of equations for $U$ and $V$:
\bea
\label{6}
\cL(t, D_t)  U(t) &=& \big(D_t+A(t)\big)U(t)={\mathcal P}f \\
\label{7}
\cL(t, D_t) V(t) &=& \big(D_t+A(t)\big)V(t)={\mathcal Q}f .
\eea
Using  representations \eqref{L1}, \eqref{10a} together with Theorem \ref{Tmain}(iii) and writing
\begin{equation*}
{\mathcal P}f(t)=\sum_{k=1}^N\sum_{j=1}^{J_k}\sum_{m=0}^{m_{k,j}-1}
f^k_{j,m  }(t)\varphi^k_{j, {m_{k,j}-1}-m}(t),
\end{equation*}
where
$$
f_{j,m  }^k (t) = \big( f(t) , \psi_{j,m  }^k(t) \big),
$$
we can present (\ref{6}) as a system of first order differential equations
\begin{equation}\label{j9a}
(D_t-\lambda_k)u^k_{j,m}(t)+u^k_{j,m-1}=f^k_{j,m},
\end{equation}
Here, $k=1,\ldots,N$, $j=1, \ldots, J_k$, $m=0, \ldots, m_{k,j}-1$,
and  $u^k_{j,m}$ is given by (\ref{KK1}) and  we assume that $u^k_{j,m-1}= 0 $ if $m=0$.

The equation (\ref{7}) concerns  the "remainder" term: here we have removed the spectrum
$\lambda_1,\ldots,\lambda_N$ from the operator using the
projector, and hence it has better estimates. This  property of the
equation (\ref{7}) is  contained in the following assertion.

\begin{theorem}\label{T9ja}
Let
$f \in  L_{\rm loc}^2(\bbR;Y)$ and
\begin{equation}\label{j3f}
\int\limits_{\bbR} \mu(t)\Vert f ; \cY(t,t+1) \Vert  dt<\infty.
\end{equation}
Then, the equation
\beas
\cL(t, D_t) u = f   
\eeas
has a solution $u = U +V \in \cX_{\rm loc}$ such that $U$ is a solution of
\ef{6} and  $V$ is a solution of \ef{7}  satisfying  the estimate
\bea
\Vert V ; \cX(\tau,\tau+1) \Vert \leq C \int\limits_{\bbR}\mu(t-\tau)\Vert {\mathcal Q}f ; \cY (t,t+1) \Vert  dt  \label{57}
\eea
for all $\tau \in \bbR$.

Let $f$ satisfy \eqref{j3f} and ${\mathcal Q}f=0$. If the bounds
\eqref{j7b} with some constant $C$ hold for $u$, then $V=0$.
\end{theorem}

Proof. We need to prove here the existence of $V$ satisfying (\ref{57}) as
well as the last uniqueness statement. We start by the uniqueness. Let
${\mathcal Q}f=0$ and assume $u$  satisfies \eqref{j7b}. Let the
coefficients of ${\mathcal P} u $  in \eqref{L1}   satisfy \eqref{j9a}. Then
$\cP u$ is a solution to (\ref{1k}) and by analysing  solutions of the
ODEs (\ref{j9a}) we conclude that this solution
also satisfies (\ref{j7b}) possibly with a slightly larger $\beta_1$
and smaller $\beta_2$. Then, according to the above uniqueness result
from Corollary \ref{Cor1a}  we have $u= {\mathcal P}u$ and hence $V={\mathcal Q}u=0$.

Let us turn to the existence. Let first $f$ have a compact support and let
$g={\mathcal Q}f$. Applying Theorem \ref{Tj6a} to the equation ${\mathcal L}V=g$, we get the estimates
\begin{equation}\label{j7c}
\int\limits_{-\infty}^\infty e^{-\beta_k (t-\tau)}\Vert V(t);X\Vert ^2dt\leq C\int\limits_{-\infty}^\infty e^{-\beta_k (t-\tau)}\Vert g(t);Y\Vert ^2dt,
\end{equation}
for $k=1,2$ and all $\tau \in \bbR$. Finally, assume $f$ is as in \ef{j3f} and
write $f=\sum_{j=-\infty}^\infty f_j$, where
$$
f_j(t)=f\ \mbox{on $(j,j+1)$ and}\ f_j=0\ \mbox{for $t$ outside $(j,j+1)$},\ j=0,\pm 1,\pm 2,\ldots.
$$
Using estimate (\ref{j7c}) for the function $V_j$ (corresponding to $f_j$ and $g_j={\mathcal Q}f_j$) we get
$$
\Big(\int\limits_\tau^{\tau+1}\Vert V_j(t);X\Vert ^2dt\Big)^{1/2}\leq C\mu(j-\tau)\Big(\int\limits_j^{j+1}\Vert {\mathcal Q}f(t);Y\Vert ^2dt\Big)^{1/2}.
$$
Summing up these relations, we obtain for all $\tau \in \bbR$
$$
\Big(\int\limits_\tau^{\tau+1}\Vert V(t);X\Vert ^2dt\Big)^{1/2}\leq C\int\limits_{\bbR}\mu(t-\tau)\Vert {\mathcal Q}f(t);Y\Vert dt,
$$
which is the same as (\ref{57}). \ \ $\Box$

\BER
We note that the function $\mu(t)$ is the Green function of the second order operator
$-(D_t-\beta_1)(D_t-\beta_2)$ up to a positive constant factor. So the estimate (\ref{57}) is similar to the representation of the solutions to the equation
$$
-(D_t-\beta_1)(D_t-\beta_2)u(t)=f(t)
$$
through the Green function and the right-hand side.

Estimates (\ref{j3f}) and (\ref{57}) imply estimate (\ref{j7b}) possibly with some slightly larger  (smaller) $\beta_1$ ($\beta_2$) and hence uniqueness and existence parts in Theorem \ref{T9ja} are in agreement.
\ENR

\section{Reduction of a dynamical system with periodic coefficients in Hilbert
space}
\label{sec7}

\subsection{Formulation of the problem}
In this section we apply Theorem \ref{T9ja}, in a slightly modified form,
to obtain a center manifold reduction for a non-linear ODE with periodic
coefficients in a Hilbert space. The proofs are mainly quite straightforward
modifications of existing result, hence, we only sketch many of them.

To formulate  the problem we let the Hilbert spaces $X$ and $Y$ and the
operators $A(t)$ and $\cA(\lambda)$ be as in Section \ref{sec2}. Then, consider
the following non-linear ODE in the space $X$,
\begin{equation}\label{N1}
{\mathcal L}(t,D_t)u(t)=f(t,\mu,u(t)),
\end{equation}
where
\bea
f : \bbR\times{\mathcal M}\times \tilde{X} \rightarrow Y   \label{N1a}
\eea
and ${\mathcal M}$ is the set of parameters, which is a neighborhood of a point
$\mu_0\in\bbR^d$, while $\tilde{X}$ is a neighborhood of zero in $X$ and
$u \in \cX_{\rm loc}$ with $u(t) \in \tilde X$ for all $t$.
Moreover we assume that the function $f$ of \eqref{N1a}
belongs to $C^{k+1}_{BU}(\bbR\times{\mathcal M}\times \tilde{X};Y)$
(the space of functions $\bbR\times{\mathcal M}\times \tilde{X}
\to Y$, which, together with all $t$-, $\mu$- and $u$-derivatives
up to order  $k+1$, $k\geq 1$, are uniformly continuous and bounded; by the $u$-derivative
we mean the standard Fr\'echet-derivative in Banach-spaces) and that
\bea
f(t,\mu_0,0)=0,\;\;\partial_u f(t,\mu_0,0)=0\;\;\mbox{for  $t\in\bbR$}.
\label{N20}
\eea
Observe that according to these equalities, $f$ is in a sense at least of the
second degree with respect  to $u$, a property which will be essential to
the fixed point argument in the proof of the main result.
We  assume also that the periodic operator function $A(t)$ belongs to $C^{k+1}_{\rm per}(\bbR;L(X,Y))$ (the space of $k+1$ times continuously
differentiable functions $\bbR \to L(X,Y)$ which are periodic in $t$).

In this section and in what follows we denote by $\lambda_1,\ldots, \lambda_N$ the eigenvalues of ${\mathcal A}(\lambda)$ on the interval $\delta_0$,
\eqref{4}, and for their multiplicities, partial multiplicities, eigenvalues and generalized eigenvalues we  keep the  notation of \ref{sec3}.
Introduce the set of indices
$$
\Theta=\{ (k,j,m)\,:\,k=1,\ldots,N;j=1,\ldots,J_k;m=0,\ldots,m_{k,j} -1 \},
$$
and the total algebraic multiplicity of all eigenvalues on the interval $\delta_0$
$$
M=\sum_{k=1}^N\sum_{j=1}^{J_k}m_{k,j}.
$$
Following  Sections \ref{sec4}, \ref{sec5} we represent the function $u$ of
\eqref{N1} as 
\begin{equation}\label{N2}
u=U+V,\;\;,\;
\end{equation}
where
\begin{equation}\label{N5}
U(t)={\mathcal P}(t)u(t)=\sum_{(k,j,m)\in\Theta}
u^k_{j,m  }(t)\varphi^k_{j, m_{k,j}-1-m}(t),
\end{equation}
with
\begin{equation}\label{N6}
u_{j,m  }^k (t) = \big( u(t) , \psi_{j,m  }^k(t) \big)
\end{equation}
and $V=u-U={\mathcal Q}u$,
${\mathcal Q}={\mathcal I}-{\mathcal P}$. Due to additional smoothness of $A(t)$ we have $\varphi^k_{j,m}\in C^k_{\rm per}(\bbR;X)$ and  $\psi^k_{j,m}\in C^k_{\rm per}(\bbR;Y)$.
Now equation (\ref{N1}) can be transformed into the system
\bea
\label{N3}
{\mathcal L}(t,D_t)U(t)&=&{\mathcal P}f(t,U(t)+V(t)),
\\
\label{N4}
{\mathcal L}(t,D_t)V(t)&=&{\mathcal Q}f(t,U(t)+V(t)).
\eea
By Theorem \ref{Tmain},  equation (\ref{N3}) can be written as the
scalar system
\begin{equation}\label{N7}
(D_t-\lambda_k)u^k_{j,m}(t)-u^k_{j,m-1}(t)=f^k_{j,m}(t,U(t)+V(t)),
\ \ \ (k,j,m) \in \Theta,
\end{equation}
where $u^k_{j,-1}=0$ and
\begin{equation}\label{M2a}
f^k_{j,m}(t,u(t))=\big( f(t,\mu,u(t)) , \psi_{j,m  }^k(t) \big)
\end{equation}
for all $k,j,m$. We introduce the notation
\begin{equation}\label{M2b}
{\bf u}(t)=\{u^k_{j,m}(t)\}_{(k,j,m)\in\Theta}
= \{ \big( u(t) , \psi_{j,m  }^k(t) \big) \}_{(k,j,m)\in\Theta}
, \ \ \ {\bf f}=\{f^k_{j,m}\}_{(k,j,m)\in\Theta}
\end{equation}
and define the 
$M \times M$-matrix 
by
\bea
{\mathcal R}\bfu =\{\lambda_ku^k_{j,m}+u^k_{j,m-1}\}_{
(k,j,m)\in\Theta}.  \label{M2c}
\eea
Given $\tau \in \bbR$ and $\xi \in \bbC^M$ we associate to the equation
\eqref{N7} the problem
\begin{equation}\label{N10a}
(D_t-{\mathcal R}){\bf u}(t)(t)={\bf f}(t,\mu,U(t)+V(t)),\;\;{\bf u}(\tau)=\xi,
\end{equation}
where $U$ is expressed through ${\bf u}$ by (\ref{N5}). This reformulation
will  replace \eqref{N3} in the later considerations (see \eqref{N7an}) in the
general case of complex Hilbert spaces $X$ and $Y$. In the next section
we make remarks on the real case.

\subsection{Projector ${\mathcal P}$ in the case of real Hilbert spaces}
\label{sec7.2}
We consider here the case when the operator $A(t)$ is real, or real valued.
Namely, in possible applications for partial differential operators having real
coefficients, it may very well be desirable to stick to the real valued
solutions, but in \eqref{N10a}, there is no guarantee that for example the
matrix $\cR$ is a real one. Accordingly, we derive here a modified equation
with a real $\cR$.

To be precise, we assume here  that our basic Hilbert spaces $X$ and $Y$, Section \ref{sec2},
are  presented as a standard complexification of the Hilbert spaces $\widehat X $,
$\widetilde X $,  $\widehat Y$ and $\widetilde Y$, respectively, over the real scalar field $\bbR$, i.e.
\beas
 X = \widehat X  \oplus i \widetilde X   \ \ \mbox{and} \ \  Y = \widehat Y \oplus i \widetilde Y
\eeas
as vector spaces. Then, $\overline x =  \widehat x - i \widetilde x$ denotes the conjugation
of an element $x = \widehat x + i \widetilde x \in \widehat X  \oplus i
\widetilde X  =  X$ (here, $\widehat x$ and $\widetilde x $ are the real and
imaginary parts of $x$). Moreover,   $x \in X$ is real, if
$x \in \widehat X $, or equivalently $x= \overline x = \widehat x $; the
notions of the real and imaginary parts of $x$ are as usual. The operator $A(t): X \to Y$ is real,
if there holds $ A(t) = \overline{ A(t)}$, where $ \overline{ A(t)}$ is the conjugation of the operator $A(t)$ defined by
\beas
\overline{ A(t)}( \widehat x + i \widetilde x ) = \overline{ A(t)( \widehat x - i \widetilde x ) } \ \ \ \mbox{for all}
\ \widehat{x},\,\tilde{x} \in X.
\eeas
If $A(t)$ is real, there holds $A(t)(\widehat X ) \subset \widehat Y$.
For example, if $X$ is a Sobolev-space $H^1 (\Omega)$ on some domain $\Omega
\subset \bbR^d$, then the above concepts have natural interpretations in
terms of real valued functions and partial differential operators
with real coefficients etc.

We introduce the operators $\widehat{P}_k: \widehat{\cY} \to \widehat{\cY}$ as
\beas
\widehat{P}_{\lambda_k}u=\sum_{j=1}^{J_k}\sum_{m=0}^{m_{k,j}-1}(u,\psi^k_{j,m})_{\widehat{\mathcal Y}}\varphi^k_{j, m_{k,j}-1-m} .
\eeas
Let also
\bea
{\mathcal P}_{\lambda_k}(t) u(t)=\sum_{j=1}^{J_k}\sum_{m=0}^{m_{k,j}-1}(u(t),\psi^k_{j,m}(t))_{Y}\varphi^k_{j, m_{k,j}-1-m}(t),
\label{Q1a}
\eea
which is now defined for all $u\in C_{\rm loc} (\bbR;Y)$.  Clearly,
$\widehat{P}_{\lambda_k}$ is the spectral projector of the operator pencil
${\mathcal L}$ corresponding to the eigenvalue $\lambda_k$ and taking values in
the space of periodic functions. Similarly to Lemma \ref{lem2.6} one can show
that the operators ${\mathcal P}_{\lambda_k}(t)$ are projectors for every $t\in\bbR$, i.e. we have ${\mathcal P}_{\lambda_k}^2={\mathcal P}_{\lambda_k}$. These projectors also have the  properties
$$
{\mathcal P}(t)=\sum_{k=1}^N{\mathcal P}_{\lambda_k}(t)
\ \ \ \mbox{and} \ \ \
{\mathcal P}_{\lambda_k}(t){\mathcal P}_{\lambda_K}(t)= \delta_{k,K}
\ \mbox{for $k,K=1,\ldots,N$}
$$

According to (\ref{j3c}), the numbers $\lambda_k-2\pi i$ are also eigenvalues of the pencil ${\mathcal A}(\lambda)$ of the same multiplicity and partial multiplicities as $\lambda_k$, and the systems
\bea
e^{-2\pi i}\varphi^k_{j,m} \ \mbox{and} \ e^{-2\pi i}\psi^k_{j,m}
, \ \ j=1,\ldots , J_k, \ m = 0,\ldots, k_{k,j}-1,
\label{Q1b}
\eea
are canonical Jordan chains corresponding to the eigenvalue $\lambda_k-2\pi i$,
and the relation (\ref{4de}) holds for them. If the projector
$\cP_{\lambda_k-2\pi i}(t)$ is defined analogously to \eqref{Q1a} by using
\eqref{Q1b}, then one can verify that actually
\begin{equation}\label{Q1}
{\mathcal P}_{\lambda_k-2\pi i}(t)={\mathcal P}_{\lambda_k}(t).
\end{equation}

In the case a real valued operator $A(t)$ the number $\overline{\lambda_k}$ is
also the eigenvalue of ${\mathcal A}(\lambda)$ of the same multiplicity and
partial multiplicities as $\lambda_k$. There also holds the important relation
$$
{\mathcal P}_{\overline{\lambda_k}}(t)=\overline{{\mathcal P}_{\lambda_k}}(t).
$$

It will be convenient to re-numerate eigenvalues as follows. We first consider
the eigenvalues $\lambda_k\in (0,\pi i)$ and numerate them by $1,\ldots,
\sigma$. Then the interval $(\pi i, 2\pi i)$ contains the same number of
eigenvalues since the numbers $\lambda_k$, $\lambda_k-2\pi i$ and
$-\lambda_k+2\pi i$ are all eigenvalues of the same multiplicity and partial
multiplicities.  The eigenvalues in
$(\pi i ,2\pi i)$  are indexed  by the numbers $-1,\ldots,-\sigma$ so that
$$
-\lambda_s=\lambda_{-s}-2\pi i.
$$
Thus,
$$
{\mathcal P}(t)=\sum_{s=1}^\sigma \big({\mathcal P}_s(t)
+ \overline{{\mathcal P}_s}(t)\big) +\big( \epsilon_0 {\mathcal P}_0(t)\big)
+ \epsilon_{\sigma+1}{\mathcal P}_{\sigma+1}(t)).
$$
Here,  $\epsilon_0 =1$, if   $2\pi i$ is an eigenvalue, otherwise
$\epsilon_0 = 0$. Similarly, $\epsilon_{\sigma+1} =1$, if   $\pi i$ is an eigenvalue, otherwise this coefficient is 0.

If $2\pi i$ is an eigenvalue, then due to (\ref{Q1})
we have  ${\mathcal P}_{2\pi i}(t)={\mathcal P}_0(t)$; one can choose
all eigenfunctions and associated eigenfunctions $\varphi^0_{m,j}$ and
$\psi^0_{m,j}$ to be real. Hence,  the expression for the operator ${\mathcal P}_0(t)$ consists of real terms only for each $t$.
If $\lambda=\pi i$ is an eigenvalue of the operator pencil ${\mathcal A}(\lambda)$ then the corresponding projector ${\mathcal P}_{\sigma+1} (t)$ has the form
\begin{equation}\label{Q4aa}
\big ({\mathcal P}_{\sigma+1}u(t)\big )(t)=\sum_{j=1}^{J_{\sigma+1}}\sum_{m=0}^{m_{\sigma+1,j}}u^{\sigma+1}_{j,m}(t)\varphi^{\sigma+1}_{j,m_{\sigma+1,j}-1-m}(t),
\end{equation}
where
\begin{equation}\label{Q4ab}
u^{\sigma+1}_{j,m}(t)=(u(t),\psi^{\sigma+1}_{j,m}(t))_Y
\end{equation}
and
\begin{equation}\label{Q4a}
\varphi^{\sigma+1}_{j,m}\ \mbox{and}\ \psi^{\sigma+1}_{j,m}, \ \ j=1,\ldots,J_{\sigma+1} ,\ m=0,\ldots,m_{\sigma+1 ,j}-1,
\end{equation}
are canonical Jordan chains  corresponding to the eigenvalue $0$ of the
operators $D_t+A(t)$ and $-D_t+A^*(t)$ with anti-periodic conditions ($u(t+1)=-
u(t)$). The chains are subject to the biorthogonality conditions (\ref{4de}).
Clearly, the functions (\ref{Q4a}) can be chosen to be real\footnote{We note
that isomorphism between periodic and anti-periodic vector-functions is given
by multiplication by $e^{\pi it}$.}. So the projector ${\mathcal P}(t)$ is real for each $t$.

To give a more explicit representation for ${\mathcal P}(t)$ we write the real
and imaginary parts
 $$
 \varphi^k_{j,m}(t)=\widehat{\varphi}^k_{j,m}+i\widetilde{\varphi}^k_{j,m}\;\;\mbox{and}\;\; \psi^k_{j,m}(t)=\widehat{\psi}^k_{j,m}+i\widetilde{\psi}^k_{j,m}.
 $$
Then for real valued $u$ we have
\bea
{\mathcal P}(t)u(t)& = & \sum_{s=1}^\sigma\sum_{j=1}^{J_s}\sum_{m=0}^{m_{s,j}-1}
\Big( \widehat{u}^s_{j,m}(t)\widehat{\varphi}^s_{j,m_{s,j}-1-m}(t)+  \widetilde{u}^s_{j,m}(t)\widetilde{\varphi}^s_{j,m_{s,j}-1-m}(t)\Big)
\rowpl
\sum_{s=0,\sigma+1}^\sigma\sum_{j=1}^{J_s}\sum_{m=0}^{m_{s,j}-1}
\epsilon_s\widehat{u}^s_{j,m}(t)\widehat{\varphi}^s_{j,m_{s,j}-1-m}(t),
\label{Q3ba}
\eea
where
$$
\widehat{u}^s_{j,m}(t)=(u(t),\widehat{\psi}^s_{j,m}(t))_Y\;\;\mbox{and}\;\;\widetilde{u}^s_{j,m}(t)=(u(t),\widetilde{\psi}^s_{j,m}(t))_Y.
$$
By taking the real and imaginary parts in (\ref{N7}) we get the real equations
\bea
\label{Q3b}
D_t\widehat{u}^s_{j,m}+\mu_s\widetilde{u}^s_{j,m}
-\widehat{u}^s_{j,m-1} &=&\widehat{f}^s_{j,m}  ,
\\
\label{Q3c}
D_t\widetilde{u}^s_{j,m}-\mu_s\widehat{u}^s_{j,m}
-\widetilde{u}^s_{j,m-1} &=&\widetilde{f}^s_{j,m},
\eea
where $\lambda_s=i\mu_s$ and
\begin{equation}\label{Q3d}
\widehat{f}^s_{j,m}(t)=(f(t),\widehat{\psi}^s_{j,m}(t))_Y \  \mbox{and} \
\widetilde{f}^s_{j,m}(t)=(f(t),\widetilde{\psi}^s_{j,m}(t))_Y.
\end{equation}
In the case $2\pi i$, or 0, are also  eigenvalues, we get
one more equation,
\begin{equation}\label{Q3e}
D_t\widehat{u}^0_{j,m}-\widehat{u}^0_{j,m-1}=\widehat{f}^0_{j,m} ,
\end{equation}
since 
$\varphi^0_{j,m}$, $\psi^0_{j,m}$ and $f^0_{j,m}$ are real, see
above. Similarly, if $\pi i$ is an eigenvalue, then there is the
additional equation
\begin{equation}\label{Q3f}
D_t\widehat{u}^{\sigma+1}_{j,m}-\widehat{u}^{\sigma+1}_{j,m-1}=\widehat{f}^{\sigma+1}_{j,m}.
\end{equation}

In the decomposition \eqref{N2}, i.e. $u(t)={\mathcal P}(t)u(t)+V(t)$, the
projector  ${\mathcal P}(t)u(t)$ has the new representation
(\ref{Q3ba}) and $V$ satisfies the same estimate (\ref{57}) as before.
However, the matrix \eqref{M2c} can now be replaced by another one, which is
obtained from the system (\ref{Q3b})--(\ref{Q3f}) and has real entries.
We keep the old notation and still write system \eqref{N10a} as
\begin{equation}\label{Aug24}
D_t{\bf u}-{\mathcal R}{\bf u}={\bf f},\;\; {\bf u}(\tau)=\xi\in\Bbb R^M,
\end{equation}
where ${\mathcal R}$ is a  real matrix. In the real case we will use the coordinates
$$
(\widehat{\xi}^s_{j,m},\tilde{\xi}^s_{j,m},\widehat{\xi}^0_{j,m},\widehat{\xi}^{\sigma+1}_{j,m}),\;\;s=1,\ldots,\sigma,
$$
where the last two are there only when $2\pi i$ or/and $\pi i$ are eigenvalues respectively.

\subsection{Model problem}
Hereinafter, we formulate all results for the case of real Hilbert spaces. Their generalizations onto complex Hilbert spaces are quite straightforward.

We will need to consider a linearization of the system \eqref{N3}--\eqref{N4}.
To formulate this we introduce for $\gamma > 0$ the weighted Sobolev-type
space ${\bf H}^1_\gamma (\bbR)$ of vector
functions ${\bf u}\,:\,\bbR\to\bbR^M$ with finite norm
$$
\|{\bf u};{\bf H}^1_\gamma (\bbR)\|_\tau=\Big(\int_{\bbR}e^{-2\gamma |t-\tau|}(\Vert{\bf u}(t)\Vert^2+\Vert D_t{\bf u}(t)\Vert^2dt\Big)^{1/2}
$$
(here and in what follows, $\Vert \cdot \Vert$ without a Banach space
denotes the Euclidean norm of $\bbR^M$)
and the space ${\bf X}_\gamma(\bbR)$ consisting of functions
$V \in \cX_{\rm loc}$ 
satisfying ${\mathcal Q}(t)V(t)=V(t)$ for almost all   $t\in\bbR$ with finite norm
$$
\|V;{\bf X}_\gamma(\bbR)\|_\tau=\Big(\int_{\bbR}e^{-2\gamma |t-\tau|} (\Vert V(t);X\Vert^2+\Vert D_tV(t);Y\Vert^2)dt\Big)^{1/2}
$$
Given $\tau \in \bbR$ and $\xi \in \bbR^M$, we pose the following linear
problem
\bea
\label{N7a}
(D_t-{\mathcal R}){\bf u}(t)&=&{\bf F}(t), \;\;{\bf u}(\tau)=\xi
\\
\label{N4a}
{\mathcal L}(t,D_t)V(t)&=& {\bf G}(t),
\eea
for the  unknown pair $(\bfu, V) \in \bfH_\gamma^1(\bbR) \times \bfX_\gamma(\bbR)$. Here,  ${\bf F}\,:\,\bbR\rightarrow {\bbR}^M$
is a given vector function belonging to  the  weighted space ${\bf L}^2_\gamma (\bbR)$ endowed  with the norm
$$
\|{\bf F};{\bf L}^2_\gamma (\bbR)\|_\tau=\Big(\int_{\bbR}e^{-2\gamma |t-\tau|}\Vert {\bf F}(t)\Vert^2dt\Big)^{1/2},
$$
and ${\bf G} \in \cY_{\rm loc}$ belongs to
the space ${\bf Y}_\gamma(\bbR)$ consisting of functions  satisfying
${\mathcal Q}(t){\bf G}(t)={\bf G}(t)$ for almost all $t\in\bbR$ and endowed
with the norm
$$
\|{\bf G};{\bf Y}_\gamma(\bbR)\|_\tau=\Big(\int_{\bbR}e^{-2\gamma |t-\tau|} \Vert {\bf G}(t);Y\Vert^2dt\Big)^{1/2}
$$

We denote by $\beta $ a positive number such that $\lambda_1, \ldots,
\lambda_N$ are the only eigenvalues of the pencil ${\mathcal A}(\lambda)$ on the intervals $\delta_{\beta'}$ with $|\beta'|\leq \beta$.

\BEL
\label{lem7.1}
Let $\gamma\in (0,\beta]$ and $\xi\in\bbR^M$, ${\bf F}\in {\bf L}^2_\gamma (\bbR)$, ${\bf G}\in {\bf Y}_\gamma(\bbR)$. Then the linear problem
(\ref{N7a}), (\ref{N4a}) has a unique solution ${\bf u}\in H^1_\gamma (\bbR)$, $V\in {\bf X}_\gamma(\bbR)$, and this solution satisfies
\begin{equation}\label{X1}
\|{\bf u};{\bf H}^1_\gamma (\bbR)\|_\tau\leq C(\Vert \xi\Vert +\|{\bf F};{\bf L}^2_\gamma (\bbR)\|_\tau)
\end{equation}
and
\begin{equation}\label{X2}
\|V;{\bf X}_\gamma(\bbR)\|_\tau\leq C\|{\bf G};{\bf Y}_\gamma(\bbR)\|_\tau,
\end{equation}
where $C$ is independent of $\tau$. In the case of complex Hilbert spaces we must assume certainly that $\xi\in \Bbb C^M$.
\ENL

Given $\gamma$ as in the lemma, we define the  linear operator
\begin{equation}\label{X3}
K: \bfL^2_\gamma (\bbR)\times \bbR^M  \times \bfY_\gamma(\bbR)
\to  H^1_\gamma (\bbR) \times {\bf X}_\gamma(\bbR)
\end{equation}
which maps the data $(\bfF, \xi , \bfG)$ into the solution of the problem (\ref{N7a})--(\ref{N4a})  satisfying the estimates (\ref{X1}), (\ref{X2}).

Let us sketch the proof of Lemma . By a shift of the $t$-variable the lemma can be
reduced to the case $\tau=0$, which is the case we will consider here.
Equation (\ref{N7a}) can be written as
$$
{\bf u}(t)=e^{{\mathcal R}t}\xi+\int_0^t e^{{\mathcal R}(t-s)}{\bf F}(s)ds.
$$
Since the eigenvalues $\lambda_k$, $k=1,\ldots,N$, are purely imaginary we have
$$
|e^{{\mathcal R}t}|\leq C(1+|t|)^{\widehat{m}-1},\ \ \widehat{m}=\max
\{ m_{k,j} \, : \, k=1, \ldots , N, j=1, \ldots, J_k \}
$$
and similar estimate is valid for the $t$-derivative. Therefore
$$
|{\bf u}(t)|\leq C\big((1+|t|)^{\widehat{m}-1}+\int_0^t (1+|t-s|)^{\widehat{m}-1}|{\bf F}(s)|ds\big)
$$
and
$$
|{\bf u}'(t)|\leq C\big((1+|t|)^{\widehat{m}-1}+\int_0^t (1+|t-s|)^{\widehat{m}-1}|{\bf F}(s)|ds+|{\bf F}(t)|\big)
$$
which implies (\ref{X1}).

We will use the following version of  Theorem \ref{T9ja} adapted to our situation.
\begin{theorem}\label{T9ja2}
Let
$f \in  L_{\rm loc}^2(\bbR;Y)$ and
\begin{equation}\label{j33f}
\int\limits_{\bbR} e^{-\beta |t|}\Vert f ; \cY(t,t+1) \Vert  dt<\infty.
\end{equation}
Then, the equation
\beas
\cL(t, D_t) u = f   
\eeas
has a solution $u = U +V \in \cX_{\rm loc}$ such that $U$ is a solution of
\ef{6} and  $V$ is a solution of \ef{7}  satisfying  the estimate
\bea
\Vert V ; \cX(\tau,\tau+1) \Vert \leq C \int\limits_{\bbR}e^{-\beta |t-\tau|}\Vert {\mathcal Q}f ; \cY (t,t+1) \Vert  dt  \label{577}
\eea
for all $\tau \in \bbR$.

Let $f$ satisfy \eqref{j33f} and ${\mathcal Q}f=0$. If the bound
\begin{equation*} 
\Vert u;{\mathcal X}(t,t+1)\Vert \leq Ce^{\beta |t|}\;\,\;\mbox{for $t\in \bbR$ }
\end{equation*}
with some constant $C$ hold for $u$, then $V=0$.
\end{theorem}

The estimate (\ref{X2}) follows from (\ref{577}) by using the arguments in
Section 4 of \cite{KM1}; the fact the in the present case the operator
$A(t)$ is not constant in $t$ does not influence the proofs.
We need to observe that $\beta$ in the reference can be replaced by $\beta_1>\beta$ with  $\beta_1$ so small enough so that the strip $\beta<\Im\lambda\leq \beta_1$ does not contain eigenvalues of the pencil ${\mathcal A}(\lambda)$.

\subsection{Center manifold reduction}

The following theorem is an analog of Theorem 1 in \cite{M1} and its proof is
literally the same. We present below the main steps of the proof for the
convenience of the reader. As in the previous section we consider the real case and hence the matrix ${\mathcal R}$ is real valued and defined in the end of Sect.\ref{sec7.2}.

\begin{theorem}\label{TKk1a} Suppose that the assumption on $A(t)$ from {\rm Sect.\ref{sec2}} are valid and the Operator $A(t)$ is real valued.
There are neighborhoods of zero $W_1\subset \bbR^M$, $W_2\subset X$ and a neighborhood ${\mathcal M}_0$ of $\mu_0$ in ${\mathcal M}$ and a function
$$
h=h(t,\mu,\xi)\in C^k(\bbR\times {\mathcal M}_0\times W_1,W_2),
$$
such that
$$
{\mathcal Q}(t)h(t)=h(t),\;\;h(t,\mu_0,0)=0\;\;\mbox{ and $\partial_\xi h(t,\mu_0,0)=0$ for $t\in\bbR$ }.
$$
Moreover the following properties hold:

1. Every solution ${\bf u}\,:\,\bbR\rightarrow W_1$ of the reduced system
\begin{equation}\label{K12a}
D_t{\bf u}-{\mathcal R}{\bf u}={\bf f}(t,\mu,U(t)+h(t,\mu,{\bf u}(t))),
\end{equation}
where $\mu\in {\mathcal M}_0$ and $U={\mathcal P}u$ is given by {\rm (\ref{Q3ba})},
gives a solution $u=U+h(t,\mu,{\bf u}(t))$ of the whole system {\rm (\ref{N1})}.

2. Every small bounded solution $u\,:\,\bbR\rightarrow X$ of {\rm (\ref{N1})}
satisfying ${\bf u}(t)\in W_1$  and $V(t)\in W_2$ for all $t\in\bbR$ yields
a solution ${\bf u}$  of the reduced system {\rm (\ref{K12a})}. In this case
$V(t)=h(t,\mu,{\bf u}(t))$.
\end{theorem}

The first step of the proof is to study the solvability of the non-linear problem
\bea
\label{N7an}
(D_t-{\mathcal R}){\bf u}(t)&=&{\bf f}(t,\mu, U(t)+V(t)),
\ \ {\bf u}(\tau)=\xi
\\
\label{N4an}
{\mathcal L}(t,D_t)V(t)&=&{\mathcal Q}(t)f(t,\mu, U(t)+V(t))
\eea
for small $\Vert \xi\Vert $. Since we are interested in small solutions it is
convenient to introduce new right-hand sides in (\ref{N7an}), (\ref{N4an})
which however coincide with the original ones for small $u$. Thus, let
\bea
\bff_\varepsilon (t,\mu, u)&=&\bff(t,\mu,u)\chi\Big(\frac{\Vert {\bf u}\Vert }{\varepsilon}\Big)\chi\Big(\frac{\Vert V\Vert_X}{\varepsilon}\Big),
\nonumber \\
f_\varepsilon (t,\mu, u)&=&f(t,\mu,u)\chi\Big(\frac{\Vert {\bf u}\Vert }{\varepsilon}\Big)\chi\Big(\frac{\Vert V\Vert_X}{\varepsilon}\Big),
\label{M1a}
\eea
where $\chi\in C^{k+1}([0,\infty))$ is a cut-off function such that
$\chi(t)=1$ for $t\leq 1$ and $\chi(t)=0$ for $t>2$ and $\varepsilon\in (0,
\varepsilon_0)$ with $\varepsilon_0$ satisfying
\bea
\{U+V\,:\,\Vert {\bf u}\Vert \leq 2\varepsilon_0,\;\Vert V\Vert_X\leq 2\varepsilon_0\}\in\widetilde{X}. \label{M1b}
\eea
Now instead of (\ref{N7an}), (\ref{N4an}) we consider the problem
\bea
\label{N7anX}
(D_t-{\mathcal R}){\bf u}(t) &=&{\bf f}_\varepsilon(t,\mu, U(t)+V(t)), \;\;{\bf u}(\tau)=\xi,
\\
\label{N4anX}
{\mathcal L}(t,D_t)V(t) &=& {\mathcal Q}(t)f_\varepsilon(t,\mu, U(t)+V(t)),
\eea
where $f_\varepsilon$ and ${\bf f}_\varepsilon$ are defined by (\ref{M1a}).
Using the operator $K$ introduced in \eqref{X3}, we can write problem (\ref{N7anX}), (\ref{N4anX}) as a fixed point equation
$$
({\bf u}(t),V(t))=K\big({\bf f}_\varepsilon(\cdot,\mu, U(\cdot)+V(\cdot)),\xi,{\mathcal Q}(\cdot)f_\varepsilon(\cdot,\mu, U(\cdot)+V(\cdot))\big)(t)
$$
Similarly to \cite{M1} and using \eqref{N20} among other things, one can show
that the operator $K(\cdot, \xi , \cdot)$ is a contraction in a small ball of the Hilbert-space $H_\gamma^1(\bbR) \times \bfX_\gamma(\bbR)$
determined by the choices \eqref{M1a}, \eqref{M1b}. The problem thus has a unique small solution for small $\varepsilon$. The function $h$ is defined by
\begin{equation}\label{M3a}
h(\tau,\mu,\xi)=V(\tau,\mu,\xi)(\tau).
\end{equation}
The same arguments as in \cite{M1} show that $h$ has the required smoothness
with respect to $t$, $\mu$ and $\xi$ and that the properties 1 and 2 hold.

\bigskip
The following remarks can be verified by quite straightforward
arguments  using the definition of the function $h$ given in the proof of the above theorem.

\BER
\label{Rem1a} If $f$ is periodic in $t$, then the function $h$ is also periodic in $t$.
\ENR

\BER
\label{Rem1} If additionally
$$
f(t,\mu,0)=0\;\;\mbox{for $\mu\in {\mathcal M}$, $t\in\bbR$},
$$
then
$$
h(t,\mu,0)=0\;\;\mbox{for $\mu\in {\mathcal M}_0$, $t\in\bbR$.}
$$
\ENR
\BER
\label{Rem2}
Let the operator $A(t)$ (and hence ${\mathcal L}$) depend also on $\mu\in{\mathcal M}$ of class $C^{k+1}_{\rm per}(\bbR\times{\mathcal M};L(X,Y))$.
So $A=A(t,\mu)$ and ${\mathcal L}(t,\mu,D_t)$. We assume also that the structure of the orthogonal projector ${\mathcal P}={\mathcal P}(t,\mu)$ is the same for $\mu\in{\mathcal M}$ and it is given by (\ref{N5}) where the functions $\varphi^k_{jm}=\varphi^k_{jm}(t,\mu)$ and $\psi^k_{jm}=\psi^k_{jm}(t,\mu)$ are smoothly depend on $\mu$ also (of class $C^{k+1}_{\rm per}(\bbR\times{\mathcal M};X)$ and $C^{k+1}_{\rm per}(\bbR\times{\mathcal M};Y)$ respectively. If
$$
f(t,\mu,0)=0\;\;\mbox{and}\;\;\partial_uf(t,\mu,0)=0\;\;\mbox{for $\mu\in {\mathcal M}$, $t\in\bbR$},
$$
then
$$
h(t,\mu,0)=0\;\;\mbox{ and $\partial_\xi h(t,\mu,0)=0$ for $\mu\in {\mathcal M}_0$ and $t\in\bbR$. }
$$
\ENR


Acknowledgement. V.\,K. was supported by the Swedish Research Council (VR),
2017-03837. J.\,T. was supported by a research grant from the Faculty of
Science of the University of Helsinki.

\end{document}